\documentclass[twocolumn]{article}
\usepackage{floatflt}
\usepackage{graphicx}
\usepackage{amsmath,latexsym}

\begin{document}
%\Large

%\baselineskip=20truept
\author{Aleksandr Koshelev \footnote{Volga branch of Moscow automobile-road institute. Email: al.koshelev@gmail.com}}
 
%\address{Volga branch of Moscow automobile-road institute.} 
%\email{freelance26@gmail.com}
%\centerline{Volga.}

\title{\textbf{Time Optimal Return of a Dynamic Object  \footnote{This paper based on my PhD thesis which I successively defend in Institute for problems in mechanics named after A.Yu. Ishlinskii Russian Academy of sciences.}}}
\date{}
\maketitle

%\begin{abstract}
\textbf{\textsl{Abstract}\textemdash We solve the problem concerning a time optimal return of a particle with a prescribed velocity to the origin by applying a magnitude-bounded force. The equations of controlled motion are derived and explicitly integrated, and the optimal open-loop control and optimal time are analyzed depending on the parameters of the problem. The qualitative behavior of the solution is established, and the solution is compared with other regimes of motion. The results are of interest for control theory and its applications to controlled flight mechanics.}
%\end{abstract}

%\begin{abstract}
%We solve the problem concerning a time optimal return of a particle with a prescribed velocity to the origin by applying a magnitude-bounded force. The equations of controlled motion are derived and explicitly integrated, and the optimal open-loop control and optimal time are analyzed depending on the parameters of the problem. The qualitative behavior of the solution is established, and the solution is compared with other regimes of motion. The results are of interest for control theory and its applications to controlled flight mechanics.
%\end{abstract}

%\text{ Time optimal control, point mechanics }


\begin{center}
\textsc{I. Introduction}
\end{center}

In the case of general initial and final conditions (synthesis), the solution to a multidimensional problem faces considerable computational difficulties. Although the equations of controlled motion can be analytically integrated in explicit form, the maximum-principle boundary value problem can be solved numerically only for modified special formulations in which the final velocity is ignored [1--3]. A numerical--analytical study with fixed magnitudes of velocity was performed for equal values of the initial and final velocities [4].

A theoretically and practically interesting case is that where the initial and final positions coincide while the corresponding velocities differ. It leads to planar (two-dimensional) motion for nonparallel velocity vectors. This plane of controlled optimal motion is determined by the indicated vectors. In the special case of parallel vectors, the problem degenerates into the one-dimensional case, for which an analytical solution can easily be constructed. In the case of two-dimensional motion, the maximum principle can be used to determine the structure of an optimal control and to completely integrate the equations of motion. The unknown parameters and the optimal time are numerically determined by Newton's method from the maximum-principle boundary conditions an arbitrary magnitude and direction of the initial velocity. Control regimes with coinciding initial and final positions (or velocities) can be used to construct quasi-optimal global motion.

\begin{center}
\textsc{II. Formulation of the problem}
\end{center}

The control system under consideration, the terminal conditions, and the functional are governed by the equations
\begin{equation}\label{E:cong}    
\dot{x}=v,       \qquad           x(0)=x_{0},  \qquad      x(t_{f})=0
\end{equation}
\[    \dot{ v}=u,          \qquad           v(0)=v_{0},  \qquad    v(t_{f})=v_f=(-1,0) \]
\[    \vert u \vert \leq 1 ,       \qquad        t_f \to \min_{u}  \qquad  x, v \in R^2\]

The more general case of arbitrary values specified for the constant mass $m$, the control force $u$ constrained by $|u| \leq u_0$, the terminal velocity $v_f$, the coinciding initial $x_0$ and final $x_f$ positions, and the initial time $t_0$ can be reduced to problem (1) by using suitable substitutions. System (1) is rotationally symmetric. In the general case where the vectors $x$ and $v$ have a dimension $n \geq 2$, the optimal control regime is equivalent to two-dimensional motion ([4,5].

The control problem has a solution in the class of piecewise continuous functions. Necessary optimality conditions of the maximum principle type [1] are applicable to this problem.

The solution of time optimal control problem (1) is reduced to solving a two-point boundary value problem for a control function of the form
\begin{equation} \label{E:cong2}   u^* = q \vert q \vert^{-1}, \qquad p(t) = p^0 = \text{const},\end{equation}
\[\qquad q(t) = -p^0 t + q^0, \qquad q^0 = \text{const}, \]
where $p$ and $q$ are the variables conjugate to $x$ and $v$ (momenta).

% \[ x(t) = x_0 + v_0 t + \int_{0}^{t} \frac{ (t-\tau) Q(t) } { R(t)} dt = x_0 + v_0 t - \frac{t}{\rho ^{2}}[-\xi R(0)+(\rho \eta - \sigma \xi)
%V(0)]+ \]
%\[ +\frac{\xi}{2 \rho ^ 3}[(\rho \tau + 3 \sigma) R(\tau) + (3 \sigma ^2 - 1)V(\tau) ]_{0}^{t} - \frac{\eta} { \rho ^2} [ \sigma V(\tau ) + R(\tau)]_{0}^{t}\]

Substituting the control function $u^*(t)$ into the equations of motion and performing single and repeated integrations, we obtain representations for the velocity $v(t)$ and the position vector $x(t)$, respectively [4, 5].

As a result, we find explicit representations of the desired phase variables $x(t)$ and $v(t)$ in terms of algebraic and logarithmic functions. They have the form of "linear" expressions in terms of vectors $\xi$ and $\eta$:
\begin{equation}\label{E:cong3}
v(t)= v_0 + \int_0^t u^* (\tau) d \tau = v_0+V_\xi (t) \xi+V_\eta (t) \eta
\end{equation}
\[ x(t)=x_0+v_0 t+X_\xi (t) \xi+X_\eta (t) \eta \]
\[ u^* = \frac{Q}{R}, \qquad Q = -\xi t + \eta \]
\[\qquad R = \vert Q \vert = (\rho^2 t^2 - 2 \sigma \rho t + 1)^{1/2} \]

Here, $\xi$ and $\eta$ are the two-dimensional vectors $p_0$ and $q_0$ normalized by $|q_0| > 0$ and $\eta$ is the unit vector. Thus, $\rho = \vert \xi \vert$, $\rho \sigma = (\xi,\eta)$, and the $2 \times 2$ matrix $(V_{\xi, \eta}, X_{\xi, \eta})$ is nonsingular in the generic case. The scalar functions $V_{\xi, \eta}$ and $X_{\xi, \eta}$ depend on time $t$ and the parameters $\rho$ and $\sigma$ in a rather intricate manner:

\begin{equation}\label{E:cong4} V_\xi = - \rho^{-2}(\sigma V+R) \vert_{0}^{t}, \qquad  V_\eta = \rho^{-1} V \vert_{0}^{t} 
\end{equation}
\[ X_{\xi} = \frac{1}{2 \rho ^3} [(- \rho \tau +3 \sigma ) R +(-2 \rho \sigma \tau + 3\sigma ^2 - 1)V]^t_0 \]
\[+ \frac{t}{\rho^2} (1+\sigma V (0))\]
\[ X_{\eta} = \rho^{-2} [ -R+ (\rho \tau - \sigma) V]^t_0 - \frac{t}{\rho} V(0)\]

The quantities to be determined are the unknown vectors $\xi$ and $\eta$ and the optimal time $t_f$.

\begin{center}
\textsc{III. Solution of the problem}
\end{center}

The boundary value problem is reduced to a system of four algebraic equations for $\xi$, $\eta$, and $t_f$ obtained from (3) and (4) taking into account boundary conditions (1). Instead of $\xi$, we introduce the vector $\zeta = \xi t_f$, where $|\zeta| = \mu = \rho t_f$. This makes it possible to separate the unknown $t_f$ in the equations. The system of two vector equations for $\zeta$, $\eta$ and $t_f$ is transformed into
\begin{equation}\label{E:cong5} 
v_f = t_{f}(a_{\zeta} \zeta + a_ \eta \eta ) 
\end{equation}
\[ v_0 - v_f = t_{f}(b_{\zeta} \zeta + b_ \eta \eta ) \]

The scalar coefficients $a_{\xi, \eta}$ and $b_{\xi, \eta}$ in (5) are functions of only $\mu$ and $\sigma$:
\begin{equation}\label{E:cong6}  
a = (\mu^{2}-2 \mu \sigma +1)^{1/2} -1,  
\end{equation}
\[ b= \textrm{ln} (( \mu - \sigma + (\mu^{2}-2 \mu \sigma +1)^{1/2}) / (1-\sigma))\]
\[ a_{\zeta}=-(1/2 \mu^{-3})[( \mu+3 \sigma) a+\mu +(3 \sigma^{2}-1)b] \]
\[ a_{\eta}=b_{\zeta}=\mu^{-2}(a+ \sigma b), \qquad b_{\eta} = -b/ \mu\]

By reducing (5) to three scalar combinations and taking into account (6), we obtain a system of three transcendental equations for $t_f$, $\mu$, and $\sigma$:
(7)

\begin{equation}\label{E:cong7} 
t_f = (f^2_x (\mu, \sigma)))^{-1/2} 
\end{equation}
\[ v^2_0 = f^2_0 (\mu, \sigma) (f^2_x (\mu, \sigma))^{-1} \]
\[ f^2_x (\mu, \sigma) \vert v_0 \vert c = - f^2_1 (\mu, \sigma), \qquad c = \cos(v_0, v_f) \]
The functions $f_x$, $f_0$, and $f_1$ depend only on $\mu \geq 0$ and $|\sigma| \leq 1$. They are obtained by taking the scalar product and have the form

\begin{equation}\label{E:cong8}
f^2_x(\mu, \sigma) = a^2_{\zeta} \mu^2 +2a_{\zeta} a_{\eta} \mu \sigma + a^2_{\eta}
\end{equation}
\[ f^2_1(\mu, \sigma) = \mu^2 a_{\zeta}(a_{\zeta} + b_{\zeta})  + \mu \sigma (a_{\zeta}(a_{\eta}+b_{\eta}) \]
\[  + a_{\eta} (a_{\eta}+b_{\eta}))  + a_{\eta} (a_{\eta} + b_{\eta})\]
\[ f^2_0(\mu, \sigma) = (a_{\zeta} + b_{\zeta})^2 \mu^2 +2 \mu \sigma (a_{\zeta} + b_{\zeta}) \]
\[ \times (a_{\eta + }b_{\eta})  + (a_{\eta} + b_{\eta})^2\]

The second and third equations in (7), combined with (8), give the governing system of equations for $\mu$ and $\sigma$. The unknown $t_f$ is then determined from the first expression. Note that a similar system was solved by a rather complicated separation of variables in [5]. That approach is also applicable to the case under consideration, but Newton's method was found to be preferable.

\begin{center}
\textsc{IV. Optimal trajectories}
\end{center}

We present the main results on the construction of optimal trajectories and examine their basic properties (Fig. 1, 2). It was established by contradiction that the optimal trajectory lies in the angle between the initial velocity and the positive abscissa axis. Moreover, the properties of $u^*(t)$ imply that the angle between the position vector $x(t)$ and the positive abscissa axis varies monotonically with time (see below).

\begin{figure}[t] 
\centering
\includegraphics[viewport=170 520 430 750,height = 5.8cm, width=7cm]{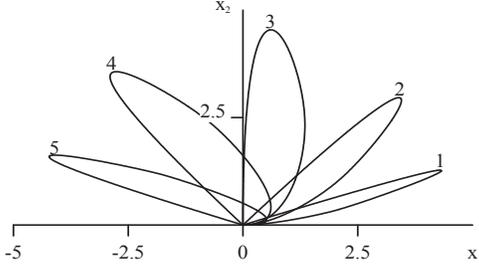}

\caption{Family of optimal trajectories for $\vert v_0 \vert = 3$ and various initial velocity directions.}
$$
\begin{array}{cccccccccc}
N & 1 & 2 & 3 &4 & 5 \\
\alpha & 0.1 \pi & 0.25 \pi & 0.5 \pi & 0.75 \pi  & 0.9 \pi  \\
\end{array}
$$
\end{figure}
\begin{figure}[t]
\centering
\includegraphics[viewport=170 520 430 750,height = 5.8cm, width=7cm]{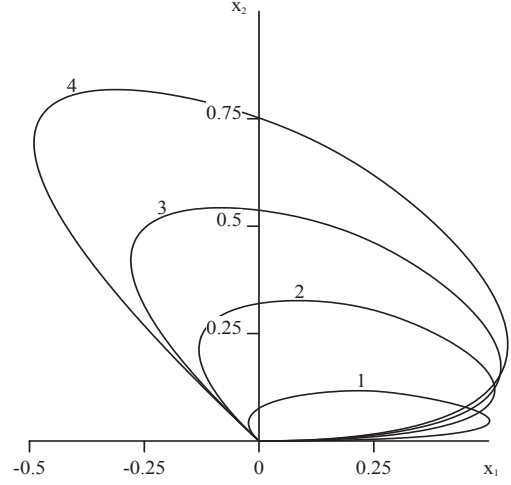}

\caption{Family of optimal trajectories for $\alpha = 0.75 \pi$ and various initial velocity magnitudes.}
$$
\begin{array}{cccccccccc}
N & 1 & 2 & 3 & 4  \\
\vert v_0 \vert & 0.3 & 0.7 & 1 & 1.3  \\
\end{array}
$$
\end{figure}

%\begin{theorem}

%Optimal control.
%\end{theorem}

After the substitution $v'_0 = v_0 \vert v_0 \vert^{-2}$, the optimal time becomes $t'_f = t_f \vert v_0 \vert^{-1}$. The parameters are transformed as follows: $\mu' = \mu$, $\sigma' = \sigma$, $\xi' = \xi \vert v_0 \vert$, and $\eta' = \eta$. When $|v_0| > 1$, the solution to problem (1) can be reduced to the case $0 < |v'_0| \leq 1$. If $|v'_0| < \vert v_0 \vert$, then control $u'(t)$ at the initial point corresponds to the control $u(t)$ at the terminal point. When $ |v_0| = 1$, the velocities and time do not vary: $v'_0 = v_0$ and $t'_f = t_f$.

The various segments of the optimal trajectories can be interpreted using the solution to the one-dimensional problem. The found features of the optimal control and the optimal time correspond the case of parallel $v_0$ and $v_f$ ($x$, $v$, and $u$ are scalar functions).

\begin{center}
\textsc{V. One-dimentional case}
\end{center}

The control function can be rewritten in the form
\begin{equation}\label{E:cong9} u=\frac{-p_0 \tau+q_{0}}{\vert -p_0 \tau+q_{0} \vert } \end{equation} 
Here, $\tau = t/t_f$ is normalized time. Switching occurs at $\tau* = 1/\lambda$, where $\lambda > 1$. In view of (9), the phase variables have the form
\[u(\tau) = \eta \text{sign}(- \lambda \tau + 1)\]
\[ \eta = \frac{v_f - v_0}{\vert v_f - v_0 \vert} \text{sign} (2 - \lambda), \qquad \zeta = \lambda \eta \]

Terminal conditions (1) for $x$ and $v$ yield the relations
\begin{equation}\label{E:cong10} x (\tau) = v_0 \tau t_f + \eta \frac{t^2_f}{2 \lambda^2} (\vert \lambda \tau - 1 \vert (1 - \lambda \tau) + 2 \lambda \tau - 1) \end{equation}
\[ v(\tau) = v_0 + \eta \frac{t_f}{\lambda (1 - \vert \lambda \tau - 1 \vert)}, \qquad \tau \in [0,1] \]

\[ t_f = \frac{\lambda (1 + \vert v_0 \vert^2 + 2 \vert v_0 \vert c)^{1/2}}{\vert 2 - \lambda \vert}, \qquad \eta = \frac{\lambda (v_0 +1)}{t_f (\lambda - 2)} \]
\[ \lambda^2 (1 - v_0) - 4 \lambda + 2 v_0 + 2=0 \]
\[ \lambda = \lambda^{\pm} = \frac{2 \pm \sqrt{2 + 2v^2_0}}{1 - v_0}, \qquad v_0 \ne 1 \]

According to (10), when $v_0 = 1$, the control $u$ is a constant vector; when $v_0 = -1$, we obtain two solutions to the maximum principle problem: optimal ($t_f = 0$, $u(t) = 0$) and nonoptimal ($t_f = 4$, $u(\tau) = \text{sign} (-2 \tau + 1)$). In the case of zero initial velocity, the optimal time is $t^0_f = 1 + \sqrt{2}$. The switching curve has the form $x = \mp \frac{v^2}{2}$ for $v_0  \geq -1$ and $v_0 < -1$, respectively.
We find the set of initial velocities corresponding to a terminal time $t_f$. For $\eta = \pm 1$, we have the expressions
\[ \lambda^{\pm} = \sqrt{2 t_f} (t_f \mp 2)^{-1/2} \]
\[ v^{\pm}_0 = -1 \pm t_f \mp \sqrt{2 t^2_f \mp 4 t_f} \]

Let us examine the functions $v^{\pm}_0 (t_f)$ when $t_f \geq 2$, since $\lambda^{\pm} < 1$ for $t_f < 2$. The minimum time $t_f = 2$ is reached at $v^{\pm}_0 = 1$. We have $v^-_0 (t_f) \geq 1$ for $t_f \geq 2$, $v^+_0 (t_f) \geq 0$ for $2 \leq t_f \leq t^0_f$, and $v^+_0(t_f) < 0$ for $t_f > t^0_f$. For the terminal time $t_f = t^0_f$, we have
\[ v^+_0 (t^0_f) = 0, \qquad v^* = v^-_0 (t^0_f) = -\sqrt{2} - 2 + \sqrt{10 + 8 \sqrt{2}} \]

Thus, one value of the optimal time corresponds to one velocity when $v_0 < 0$ or $v_0 > v^*$ and to two velocities in one direction when $0 \leq v_0 \leq v^*$ (see Section VI). The found properties of one-dimensional motion are used to interpret the solution to the problem in the general (two-dimensional) case.

\begin{center}
\textsc{VI. Numerical simulation}
\end{center}

Figures 1 and 2 show the families of optimal trajectories. The optimal trajectories lie in the angle between the initial velocity direction and the positive abscissa axis. The angle characterizing the direction of $x(t)$ is a monotone function of time.

The optimal control is analyzed for various initial velocities. The resulting family makes it possible to find the optimal control and to determine its features. Note the properties of the open-loop control $u^*(t)$ at the initial and terminal times of motion (Fig. 3, 4). 

We analyze these values as functions of the initial velocity direction, i.e., the angle $\alpha$ for fixed velocity magnitudes. The desired angles $\phi$ (i.e., the directions of $u$) at the initial time are in the triangle bounded by the lines $\phi = 2 \pi$, $\phi = \alpha + \pi$, and $\alpha = 0$ (Fig. 3). The object accelerates (i.e., the velocity magnitude increases) when $\phi \geq \alpha + 1.5 \pi$ and slows down when $ \alpha + \pi \leq \phi < \alpha + 1.5 \pi$. Thus, deceleration occurs when $|v_0| \in [1; +\infty)$ and $\alpha \in [0, \pi]$ and when $\alpha \in [\pi/2, \pi]$ and $|v_0| \geq 0$ at the initial time. When $\alpha \in [0, \pi/2)$ and $|v_0| \in [0, 1)$, the motion switches from acceleration to deceleration as $\alpha$ increases. According to the solution of the one-dimensional problem, when $|v_0| \in [1, +\infty)$, we have $\phi \to \pi$ as $\alpha \to 0$ and $\phi \to 2 \pi$ as $\alpha \to \pi$; when $|v_0| \in [0, 1)$, we have $\phi \to 2 \pi$ as $\alpha \to 0$ or $\alpha \to \pi$ (see Section V).

\begin{figure}[t] 
\centering
\includegraphics[viewport=170 520 430 750,height = 5.8cm, width=7cm]{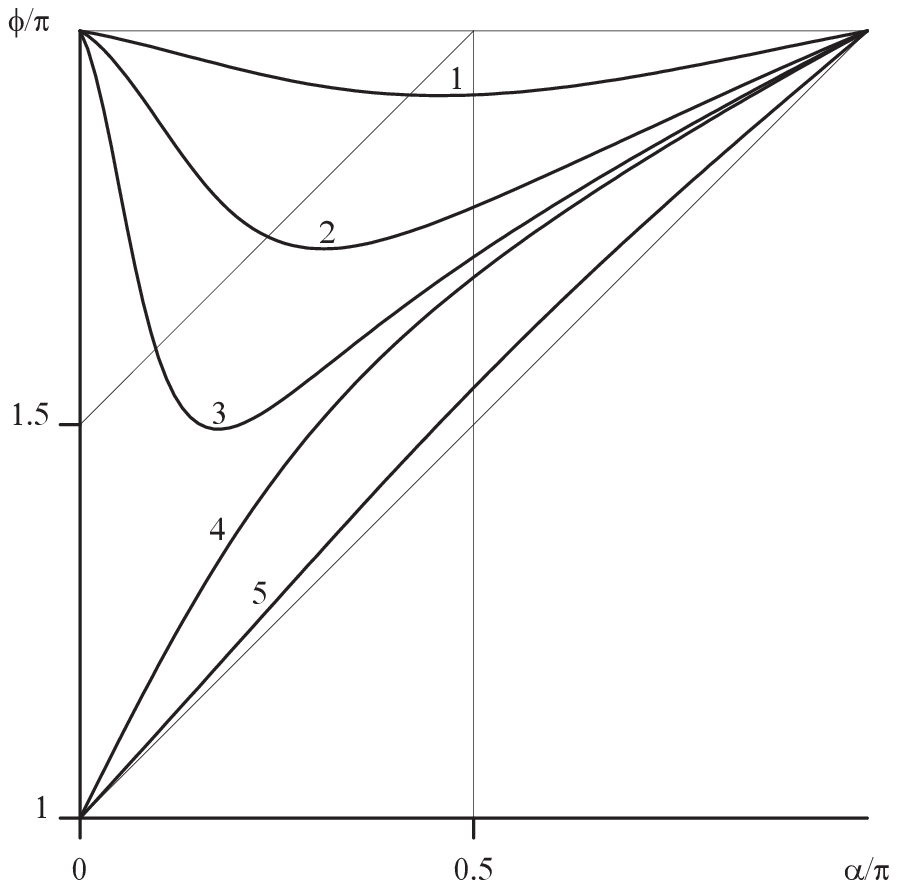}

\caption{Initial control directions $\phi$ as functions of $\alpha$ for various initial velocity magnitudes }
$$
\begin{array}{cccccccccc}
N & 1 & 2 & 3 & 4 & 5  \\
\vert v_0 \vert & 0.5 & 0.7 & 0.9 & 1 & 3  \\
\end{array}
$$

\end{figure}

\begin{figure}[t] 
\centering
\includegraphics[viewport=170 520 430 750,height = 5.8cm, width=7cm]{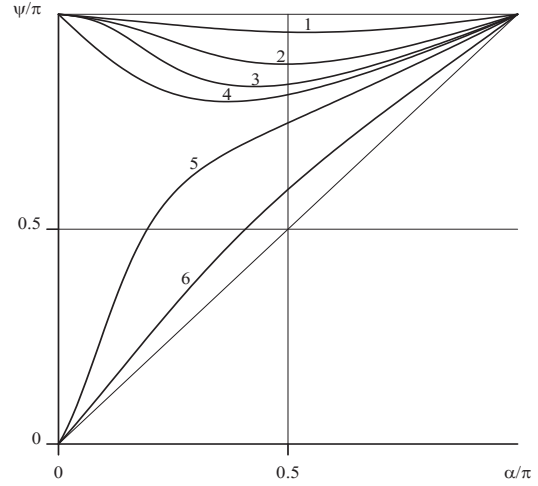}

\caption{Final control directions $\psi$ as functions of $\alpha$ for various initial velocity magnitudes }

$$
\begin{array}{cccccccccc}
N & 1 & 2 & 3 & 4 & 5 & 6 \\
\vert v_0 \vert & 0.3 & 0.7 & 0.9 & 1 & 1.3 & 3  \\
\end{array}
$$
\end{figure}

A similar diagram can be constructed for $t_f$ (Fig. 4). The range of $\alpha$ is bounded by the lines $\psi = \pi$, $\psi = \alpha$, and $\alpha = 0$. The transition between the deceleration and acceleration regimes for obtaining the required velocity $v_f$ occurs on the line $\psi = \pi/2$. When $\alpha \in [\pi/2, \pi]$ and $|v_0| > 0$ and when $|v_0| \in [0, 1]$ and $\alpha \in [0, \pi]$, the object accelerates toward the required terminal velocity. When $|v_0| \in [0, 1]$, we have $\psi \to \pi$ as $\alpha \to 0$ or $\alpha \to \pi$. When $|v_0| \in (1, +\infty)$, $\psi \to 0$ as $\alpha \to 0$ and $\psi \to \pi$ as $\alpha \to \pi$.

\begin{figure}[t] 
\centering
\includegraphics[viewport=170 520 430 750,height = 5.8cm, width=7cm]{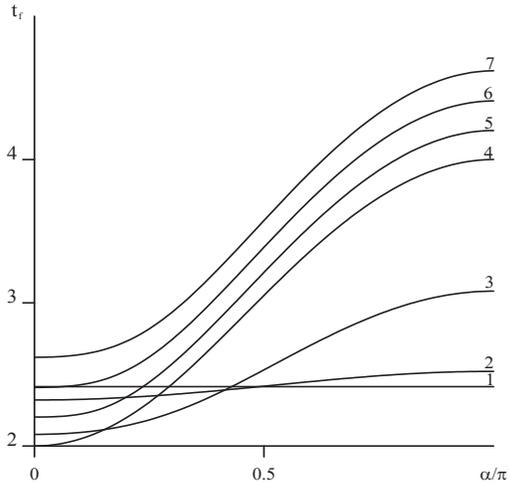}

\caption{Optimal time as a function of the initial velocity direction $\alpha \in [0, \pi]$ for various initial velocity magnitudes }
$$
\begin{array}{cccccccccc}
N & 1 & 2 & 3 & 4 & 5 & 6 & 7 \\
\vert v_0 \vert & 0 & 0.1 & 0.5 & 1 & 1.1 & 1.2 & 1.3  \\
\end{array}
$$
\end{figure}

\begin{figure}[t] 
\centering
\includegraphics[viewport=170 520 430 750,height = 5.8cm, width=7cm]{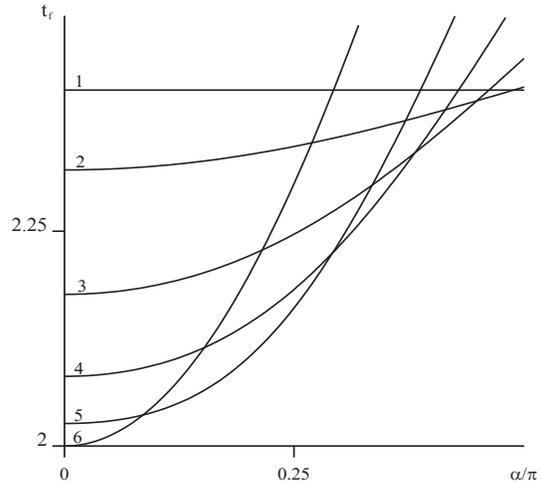}

\caption{Optimal time as a function of the initial velocity direction $\alpha \in [0, \pi/2]$ for various initial velocity magnitudes }
$$
\begin{array}{cccccccccc}
N & 1 & 2 & 3 & 4 & 5 & 6  \\
\vert v_0 \vert & 0 & 0.1 & 0.3 & 0.5 & 0.7 & 1   \\
\end{array}
$$
\end{figure}

Thus, for high initial velocities $(v_0 > 1)$, the behavior of the control at the initial stage is similar to the case of low velocities $(v_0 < 1)$ at the final stage of motion. The constraints on the direction of the control vector imply that the velocity component along the ordinate axis decreases at the initial time (Fig. 3). Moreover, there are no optimal controls under which the object accelerates at the initial stage of motion and decelerates in the neighborhood of the terminal point.

Now, we analyze the optimal time as a function of the initial velocity direction $\alpha$ (Fig. 5,6). When $\alpha = 0$ or $\pi$, the optimal time corresponds to the solution of the one-dimensional problem in Section V. Moreover, when $|v_0| = 1$, $t_f(\alpha) \to 4$ as $\alpha \to \pi$, i.e., the optimal time tends to a nonoptimal solution. In the interval $\alpha \in [0, \pi]$, the function $t_f(\alpha)$ increases monotonically. When $\alpha \in [0, \pi/2]$ and $v_0 \in [0, v*]$, there is a domain where $t_f \leq f^0_f$ (Fig. 5). In this domain, the same $t_f$ can be reached at two different velocity magnitudes. Moreover, any two curves intersect when $0 \leq |v_0| \leq 1$.

A transition occurs at the intersection point: as $\alpha$ increases, the curve with a higher velocity corresponds to a longer time; and the curve with a lower velocity, to a shorter time. The properties of the control-switching curve (see Fig. 3) and the envelope of the family of curves $t_f(\alpha)$ (Fig. 5,6) imply that, at the intersection point, deceleration is an optimal control strategy for an object with an higher velocity at the initial stage and acceleration is optimal for the control of an object with an initially lower velocity. The envelope is the control-switching curve at the initial point (Fig. 3), i.e., for fixed $\alpha \in [0, \pi/2]$. The minimum optimal time is reached when $v_0$ is orthogonal to $\eta$. For fixed values $\alpha \in [\pi/2, \pi]$, the minimum time is reached at $v_0 = 0$.

\begin{center}
\textsc{Conclusions}
\end{center}

These results make it possible to draw an analogy with the case of one-dimensional motion. The return of the object from the initial position to the same one with turning the velocity can be divided into two stages: motion from the initial position to an intermediate point $x_1$ with at the velocity $v_1$ and motion from this point to the final position. In the case of one-dimensional motion, acceleration and deceleration are possible at each stage. These regimes determine the properties of the trajectories and the control in the two-dimensional problem.

These features are represented as acceleration and deceleration regimes. In the one-dimensional case, the velocity magnitude can have three extrema and up to four basic control regimes are possible in the course of this motion. When $\alpha = 0$, $\phi = 0$, and $0 < |v_0| < 1$, we have acceleration, the reversion of the control's sign, deceleration, and again acceleration. When $\alpha = 0$, $\phi = \pi$, and $|v_0| > 1$, we have deceleration (the initial velocity magnitude decreases), then acceleration, the reversion of the control's sign, and deceleration. When $\alpha = 0$, $\phi = \pi$, and $|v_0| = 1$, the control is a constant vector implementing deceleration and acceleration. When $\alpha = 0$, we have $\phi = 0$ for all the velocities $|v_0| \geq 0$ with $|v_0| \ne 1$. Deceleration is followed by acceleration, then the sign of the control reverses, and again deceleration is followed by acceleration.

Thus, the velocity magnitude attains a minimum when $\alpha = 0$, $\phi = 0$, and $|v_0| = 1$ and a maximum and two minima when $\alpha = \pi$, $\phi = 0$, and $|v_0| \geq 0$ with $|v_0| \ne 1$. A maximum and a minimum occur in the remaining cases.

In the two-dimensional case when $|v_0| = \text{const}$ and $\alpha \in (0, \pi)$, the control leads to switching the acceleration and deceleration regimes corresponding to the boundary values of $\alpha$. This switch occurs at the initial stage of the motion if $|v_0| < 1$ and at the final stage if $|v_0| > 1$ (Fig. 3,4).

\begin{center}
\textsc{References}
\end{center}

\begin{description}
\item{[1]}
        L. S. Pontryagin, V. G. Boltyanskii, R. V. Gamkrelidze, and E. F. Mishchenko, \emph{The Mathematical Theory of Optimal Processes}             (Nauka, Moscow, 1969; Gordon  Breach, New York, 1986).
\item{[2]}
        G. Leitmann \emph{An Introduction to Optimal Control}, McGraw-Hill, NY, 1966
\item{[3]}
        A. Bryson; Ho, Yu-Chi, \emph{Applied Optimal Control : Optimization, Estimation and Control}, John Wiley \& Sons, 1975. Revised          Printing, 481 pp.
\item{[4]}
        L. D. Akulenko and A. P. Koshelev, \emph{Time-Optimal Steering of a Dynamic Object to a Given Position 
under the Equality of the Initial and Final Velocities}, Journal of Computer and System Sciences International, p. 921 - 928, Vol. 42, No. 6, 2003

\item{[5]}
        L. D. Akulenko, A. P. Koshelev, \emph{Time-optimal steering of a point mass to a specified position with the required velocity},  Journal         of applied mathematics and mechanics, Pages 200-207, Volume 71, Issue 2, 2007
\end{description}

\end{document}